\newcommand{\abstand}{\vspace{1em}}

\newtheorem{theo}{Theorem}
\newtheorem{lemma}{Lemma}
\newtheorem{remark}{Remark}

\newcommand{\proof}{{\em Proof. }}
\newcommand{\qed}{$\hfill\Box$}

\newcommand{\PG}[2]{\mbox{$\mbox{{\rm PG}}(#1,#2)$}}

\newcommand{\GF}[1]{\mbox{$\mbox{{\rm GF}}(#1)$}}

\newcommand{\abb}[3]{\mbox{$#1\,:\,#2\rightarrow#3$}}
\newcommand{\Abb}[5]{\mbox{$#1\,:\,#2\rightarrow#3,\;#4\mapsto #5$}}
\newcommand{\im}{\mbox{\rm im\,}}
\newcommand{\id}{\mbox{\rm id}}

\newcommand{\spn}{\mbox{{\rm span\,}}}
\newcommand{\Aut}{\mbox{{\rm Aut\,}}}

\newcommand{\qu}[1]{\overline{#1}}
\newcommand{\inv}{^{-1}}

\newcommand{\abf}{{\bf a}}
\newcommand{\ebf}{{\bf e}}
\newcommand{\vbf}{{\bf v}}
\newcommand{\wbf}{{\bf w}}
\newcommand{\Cbf}{{\bf C}}
\newcommand{\Ebf}{{\bf E}}
\newcommand{\Gbf}{{\bf G}}
\newcommand{\Vbf}{{\bf V}}
\newcommand{\Wbf}{{\bf W}}

\newcommand{\Bcal}{{\cal B}}
\newcommand{\Ccal}{{\cal C}}
\newcommand{\Ecal}{{\cal E}}
\newcommand{\Hcal}{{\cal H}}
\newcommand{\Kcal}{{\cal K}}
\newcommand{\Lcal}{{\cal L}}
\newcommand{\Mcal}{{\cal M}}
\newcommand{\Pcal}{{\cal P}}
\newcommand{\Qcal}{{\cal Q}}
\newcommand{\Scal}{{\cal S}}
\newcommand{\Vcal}{{\cal V}}

\documentclass[a4paper,12pt]{article}
\usepackage{amssymb}
\begin{document}

\title{The Veronese Surface in \PG53 \\and Witt's $5$--$(12,6,1)$ Design%
\thanks{Research supported by the Austrian FWF, project P12353--MAT.}
}

\author{%
Hans Havlicek\\
Abteilung f\"ur Lineare Algebra und Geometrie\\
Technische Universit\"at\\
Wiedner Hauptstra{\ss}e 8--10\\
A--1040 Wien\\
Austria\\
EMAIL: {\tt havlicek@geometrie.tuwien.ac.at}
}
\date{}
\maketitle

%
%
%
%

   \begin{abstract}
   A conic of the Veronese surface in \PG53 is a quadrangle. If one such
   quadrangle is replaced with its diagonal triangle, then one obtains a
   point model $\Kcal$ for Witt's $5$--$(12,6,1)$ design, the blocks being
   the hyperplane sections containing more than three (actually six) points
   of $\Kcal$. As such a point model is projectively unique, the present
   construction yields an easy coordinate--free approach to some results
   obtained independently by H.S.M.\ Coxeter and G.\ Pellegrino, including a
   projective representation of the Mathieu group $M_{12}$ in \PG53.
   \end{abstract}


\section{Introduction}

Throughout this paper $\Vbf$ is a $3$--dimensional vector space over
$F:=\GF3$ and $\Wbf$ denotes the symmetric tensor product $\Vbf\vee\Vbf$.
Occasionally, it will be convenient to use coordinates. We fix an ordered
basis $(\ebf_0,\ebf_1,\ebf_2)$ of $\Vbf$. It yields the ordered basis
   \begin{displaymath}
   (\ebf_0\vee\ebf_0, 2 \ebf_0\vee\ebf_1, 2 \ebf_0\vee\ebf_2,
    \ebf_1\vee\ebf_1, 2 \ebf_1\vee\ebf_2, \ebf_2\vee\ebf_2)
   \end{displaymath}
of $\Wbf$. All coordinate vectors are understood with respect to one of these
bases. The projective plane on $\Vbf$ is
$\PG23=(\Pcal(\Vbf),\Lcal(\Vbf),\in)$, where $\Pcal(\Vbf)$ and $\Lcal(\Vbf)$
denote the sets of points and lines, respectively. Likewise we have
$\PG53=(\Pcal(\Wbf),\Lcal(\Wbf),\in)$. The {\em Veronese mapping}\/ is given
by
   \begin{displaymath}
   \Abb{\varphi}{\Pcal(\Vbf)}{\Pcal(\Wbf)}{F\abf}{F(\abf\vee\abf)}
   \end{displaymath}
or, in terms of coordinates, by
   \begin{equation} \label{vero}
   F(x_0,x_1,x_2)\mapsto F(x_0^2,x_0 x_1,x_0 x_2,x_1^2,x_1 x_2,x_2^2).
   \end{equation}
The set $\im\varphi$ is the well--known
{\em Veronese surface}. See, among others, \cite[Chapter V]{bura61},
\cite{herz82}, \cite[Chapter 25]{hirs-thas91}. Recall three major properties
of the Veronese mapping:
Firstly, $\varphi$ is injective.
Secondly, the $\varphi$--image of each line
$l$ of $\PG23$ is a (non--degenerate) conic or, in other words, a planar
quadrangle in $\PG53$. The plane of this conic meets $\im\varphi$ in exactly
four points. Each conic of $\im\varphi$ arises in this way.
Thirdly, the pre--image under $\varphi$ of each hyperplane $\Hcal$ of $\PG53$
is a (possibly degenerate) quadric of $\PG23$. Each quadric of $\PG23$ arises
in this way.

If we are given a quadrangle $\Gamma$ in a projective plane of order $3$,
then its diagonal points form a triangle $\Delta$, say. On the other hand, if
$\Delta$ is a triangle in such a plane, then there are exactly four points
which are not on any side of $\Delta$. Those four points form a quadrangle,
say $\Gamma$, which in turn has $\Delta$ as its diagonal triangle
\cite[391--392]{hirs79}. This one--one correspondence between quadrangles and
triangles in a projective plane of order three is the backbone of this paper.

There is also another interpretation of this correspondence: We may consider
the quadrangle $\Gamma$ as a conic. It will be called the {\em associated
conic}\/ of the triangle $\Delta$.
The internal points of the conic $\Gamma$ comprise the triangle $\Delta$.
Moreover, $\Delta$ is a self polar triangle of $\Gamma$ \cite[Theorem
8.3.4.]{hirs79}. Finally, the sides of $\Delta$ are all the external lines of
$\Gamma$.

\section{Variations on $13-4+3=12$}

In the sequel an arbitrarily chosen line $l_\infty$ of \PG23 will be regarded
as {\em line at infinity}. Its Veronese image
$l_\infty^\varphi=:\Gamma_\infty$ is a planar quadrangle with diagonal
triangle $\Delta_\infty$, say. The plane spanned by $\Gamma_\infty$ is
denoted by $\Ecal_\infty$.

The following Theorem describes the essential construction:

   \begin{theo}\label{theo-a}
   Write $\Kcal$ for that set of points in\/ \PG53 which is obtained from the
   Veronese surface\/ $\im\varphi$ by replacing the planar quadrangle\/
$\Gamma_\infty$,
   i.e. the $\varphi$--image of the line at infinity, with its diagonal
   triangle $\Delta_\infty$. Then the following hold true:
   \begin{eqnarray}
   \label{condi}
   &d_\Hcal:=\#(\Hcal\cap\Kcal)\in\{0,3,6\}
   \mbox{ for all hyperplanes } \Hcal \mbox{ of\/ } \PG53.&\\
   \label{condii}
   &\#\Kcal=12.&
   \end{eqnarray}
   \end{theo}
\proof
The pre--image of $\Hcal$ under $\varphi$ is a quadric of \PG23, say $\Qcal$.
There are four cases \cite[140]{hirs79}.
      \begin{enumerate}
      \item
      $\Ecal_\infty\subset \Hcal$: Hence $d_\Hcal=\#\Qcal-4+3$. As
      $l_\infty\subset\Qcal$, we obtain that $\Qcal$ is the repeated line
      $l_\infty$ or a cross of lines. Thus $d_\Hcal=4-4+3=3$ or
      $d_\Hcal=7-4+3=6$.
      \item
      $\Ecal_\infty\cap \Hcal$ is an external line of $\Gamma_\infty$:
      Hence $d_\Hcal=\#\Qcal-0+2$. As $\Qcal$ is either a single affine point
      or a conic without points at infinity, we infer
      $d_\Hcal=1-0+2=3$ or $d_\Hcal=4-0+2=6$.
      \item
      $\Ecal_\infty\cap \Hcal$ is a tangent of $\Gamma_\infty$:
      A tangent carries no internal points so that $d_\Hcal=\#\Qcal-1+0$. The
      quadric $\Qcal$ is a repeated line $l$ with $l\neq\l_\infty$, or a
      cross of lines with double point at infinity, but each line other than
      $l_\infty$, or a conic touching $l_\infty$. Thus
      $d_\Hcal=4-1+0=3$ or $d_\Hcal=7-1+0=6$ or $d_\Hcal=4-1+0=3$.
      \item
      $\Ecal_\infty\cap \Hcal$ is a bisecant of $\Gamma_\infty$: A bisecant
      carries exactly one internal point, whence $d_\Hcal=\#\Qcal-2+1$. Now
      $\Qcal$ is a cross of lines with double point not at infinity, or a
      conic with two distinct points at infinity. Hence $d_\Hcal=7-2+1=6$ or
      $d_\Hcal=4-2+1=3$.
      \end{enumerate}
Finally, $\im\varphi\cap\Ecal_\infty=\Gamma_\infty$ implies
$\#\Kcal=13-4+3=12$.
\qed

   \begin{remark}
   {\em
   If $l_\infty$ is chosen to be the line $x_0=0$, then $\Delta_\infty$ can
   easily be expressed in terms of coordinates as
   \begin{equation}\label{delta}
    \{ F(0,0,0,1,0,1),\;
       F(0,0,0,2,1,1),\;
       F(0,0,0,2,2,1)\}.
   \end{equation}
   Thus, by virtue of (\ref{vero}) and (\ref{delta}), one may describe
   $\Kcal$ in terms of coordinates.
   }
   \end{remark}
Before we are going to reverse the construction of Theorem \ref{theo-a}, we
prove the following

   \begin{lemma}\label{lemma-a}
   Let $\Kcal$ be a set of points in\/ \PG53. Then (\ref{condi})
   and (\ref{condii}) together are equivalent to the conjunction of the
   following three conditions:
      \begin{eqnarray}
      \label{condv}
      &\mbox{Any\/ $5$--subset of $\Kcal$ is independent.}&\\
      \label{condvi}
      &\mbox{$\#(\Hcal\cap\Kcal)\geq 5$ implies $\#(\Hcal\cap\Kcal)=6$ for
      all hyperplanes $\Hcal$ of\/ \PG53.}&\\
      \label{condvii}
      &\#\Kcal\geq7.&
      \end{eqnarray}
   \end{lemma}
\proof {\em (\ref{condi}) and (\ref{condii}) $\Longrightarrow$ (\ref{condv})
and (\ref{condvi}) and (\ref{condvii})}:
Choose any $5$--set $\Mcal\subset\Kcal$ and $P\in\Kcal\setminus\Mcal$. At
first we are going to show that
   \begin{equation}\label{b-1}
   \dim\spn(\Mcal\cup\{P\}))\geq 4;
   \end{equation}
here ``$\dim$'' denotes the projective dimension. Assume to the contrary that
$\dim\spn(\Mcal\cup\{P\}) < 4$. Then each hyperplane of \PG53 passing through
$\Mcal\cup\{P\}$ meets $\Kcal$ in exactly six points, by (\ref{condi}). All
those hyperplanes are covering $\Kcal$, whence $\Kcal=\Mcal\cup\{P\}$, in
contradiction to (\ref{condii}).

We infer from (\ref{b-1}) that $\dim\spn\Mcal\geq 3$. This dimension cannot
equal three, since then $\Kcal$ would only have nine points, namely the five
points in $\Mcal$ plus one more point in each of the four hyperplanes through
$\Mcal$. Consequently, $\Mcal$ is independent. By (\ref{condi}) and
(\ref{condii}), conditions (\ref{condvi}) and (\ref{condvii}) follow
immediately.

{\em (\ref{condv}) and (\ref{condvi}) and (\ref{condvii}) $\Longrightarrow$
(\ref{condi}) and (\ref{condii})}:
By our assumptions, $\Kcal$ contains a basis $\Scal$ of \PG53. Each of the
six hyperplane faces of that basis contains exactly one more point of
$\Kcal$; it is in general position with respect to the remaining five. Thus
we have $\#\Kcal\geq 12$. On the other hand choose four points in $\Scal$.
Each of the four hyperplanes passing through them meets $\Kcal$ in at most
six points. Hence $\#\Kcal\leq 12$. Thus (\ref{condii}) holds true.

If we fix one $3$--set $\Delta\subset\Kcal$, then the number hyperplanes
through $\Delta$ is $13$, and the number of $2$--sets in
$\Kcal\setminus\Delta$ is $36$. By (\ref{condv}) and (\ref{condvi}), the
number of hyperplanes through $\Delta$, meeting $\Kcal$ in exactly six
points, is $36/3=12$. Hence there is a unique hyperplane $\Hcal_\Delta$, say,
with
   \begin{equation}
   \Delta=\Kcal\cap\Hcal_\Delta.
   \end{equation}

Next fix one point $P\in\Kcal$. There are $330$ $4$--subsets of
$\Kcal\setminus\{P\}$. They give rise to the $330/5=66$ hyperplanes through
$P$ meeting $\Kcal$ in six points. Likewise one finds ${11\choose 2}=55$
triangles in $\Kcal$ containing $P$. Each of those triangles yields exactly
one hyperplane through $P$ meeting $\Kcal$ in three points only. There are,
however, only $121=66+55$ hyperplanes through $P$, whence (\ref{condi})
follows.
\qed

\abstand
\noindent
Theorem \ref{theo-a} can be reversed now as follows:

   \begin{theo}\label{theo-c}
   Let $\Kcal$ be a set of points in\/ \PG53 satisfying (\ref{condi}) and
   (\ref{condii}). Suppose that $\Vcal$ is obtained from $\Kcal$ by replacing
   one triangle $\Delta\subset\Kcal$ with its associated conic\/ $\Gamma$.
   Then $\Vcal$ is projectively equivalent to the Veronese surface
   $\im\varphi$.
   \end{theo}
\proof By Lemma \ref{lemma-a}, there is a triangle $\Delta\subset \Kcal$. The
plane spanned by $\Delta$ is denoted by $\Ecal$. According to \cite[Theorem
25.3.14]{hirs-thas91} it is sufficient to verify the following conditions:
   \begin{eqnarray}
   \label{condiii}
   &c_\Hcal:=\#(\Hcal\cap\Vcal)\in\{1,4,7\}
   \mbox{ for all hyperplanes } \Hcal \mbox{ of } \PG53.&\\
   \label{condiv}
   &c_{\Hcal_0}=7 \mbox{ for some hyperplane }\Hcal_0 \mbox{ of }\PG53.&
   \end{eqnarray}

In order to establish (\ref{condiii}) choose a hyperplane $\Hcal$ and put
$d_\Hcal:=\#(\Hcal\cap\Kcal)$. There are four cases.
      \begin{enumerate}
      \item
      $\Ecal\subset \Hcal$: By (\ref{condi}),
      $c_\Hcal=d_\Hcal-3+4\in\{1,4,7\}$.
      \item
      $\Ecal\cap \Hcal$ is an external line of $\Gamma$: Thus
      $\#(\Hcal\cap\Delta)=2$ and $c_\Hcal=d_\Hcal-2+0\in\{1,4\}$.
      \item
      $\Ecal\cap \Hcal$ is a tangent of $\Gamma$: Thus
      $\#(\Hcal\cap\Delta)=0$ and $c_\Hcal=d_\Hcal-0+1\in\{1,4,7\}$.
      \item
      $\Ecal\cap \Hcal$ is a bisecant of $\Gamma$: Thus
      $\#(\Hcal\cap\Delta)=1$ and $c_\Hcal=d_\Hcal-1+2\in\{4,7\}$.
      \end{enumerate}
Two points in $\Kcal\setminus\Delta$ together with $\Delta$ generate a
hyperplane $\Hcal_0$ meeting $\Kcal$ in six distinct points by
(\ref{condv}). According to case 1, $c_{\Hcal_0}=7$.
\qed

\abstand
\noindent
All properties of the Veronese surface that are used in the following proof
can be read off, e.g., from \cite[Section 25.1]{hirs-thas91}.

   \begin{theo}\label{theo-d}
   Suppose that $\Kcal$, $\Kcal'$ are sets of points in\/ \PG53 subject to
   (\ref{condi}) and (\ref{condii}). Choose five distinct points
   $P_0,\ldots,P_4$ in $\Kcal$ and five distinct points $P_0',\ldots,P_4'$ in
   $\Kcal'$. Then there is a unique collineation $\kappa$ of\/ \PG53 with
   $\Kcal^\kappa=\Kcal'$ and $P_i^\kappa=P_i'$ for $i=0,\ldots,4$.
   \end{theo}
\proof
   Put $\Delta:=\{P_0,P_1,P_2\}$. Define $\Gamma$ and $\Vcal$ according to
   Theorem \ref{theo-c}. Write $\Ccal$ for the set of all conics
   contained in $\Vcal$. Then $(\Vcal,\Ccal,\in)$ is a projective plane of
   order 3. Moreover, the Veronese mapping $\varphi$ yields a collineation of
   \PG23 onto that projective plane. There is a unique conic in $\Vcal$
   joining $P_3$ with $P_4$. It meets $\Gamma$ in a single point, say $G_3$.
   The line spanned by $G_3$ and $P_i$ ($i=0,1,2$) is a bisecant of
   $\Gamma$, as it contains the internal point $P_i$; hence it meets
   the conic $\Gamma$ residually in a point $G_i$, say. Thus
   $\Gamma=\{G_0,\ldots,G_3\}$. The four points $\{P_3,P_4,G_0,G_1\}$ form a
   ``quadrangle'' of the projective plane $(\Vcal,\Ccal,\in)$, i.e. a set
   of four points no three of which are on a common conic  $\subset\Vcal$.

   Repeat the previous construction with $\Kcal'$ to obtain $\Delta'$ etc.
   By Theorem \ref{theo-c}, there exists a collineation $\mu$ of \PG53 with
   $\Vcal^\mu=\Vcal'$. Thus $\{P_3^\mu,P_4^\mu,G_0^\mu,G_1^\mu\}$ is a
   ``quadrangle'' of the projective plane $(\Vcal',\Ccal',\in)$. There is a
   projective collineation $\lambda'$ of $(\Vcal',\Ccal',\in)$ with
      \begin{displaymath}
      P_3^\mu\mapsto P_3', \, P_4^\mu\mapsto P_4',\,
      G_0^\mu\mapsto G_0', \, G_1^\mu\mapsto G_1'.
      \end{displaymath}
   This $\lambda'$ extends to a projective collineation $\lambda$ of \PG53.
   The product $\kappa:=\mu\lambda$ has the required properties, since
   $G_3^\kappa=G_3'$ implies $G_2^\kappa=G_2'$, so that also
      \begin{displaymath}
      P_i^\kappa=P_i' \mbox{ for } i=0,1,2.
      \end{displaymath}
   If $\qu\kappa$ is a collineation subject to the conditions of the theorem,
   then $\qu\kappa\kappa\inv$ restricts to a collineation of
   $(\Vcal,\Ccal,\in)$ fixing each point of a ``quadrangle''. Now
   $\Aut\GF3=\{\id\}$ forces $\qu\kappa\kappa\inv$ to fix $\Vcal$ pointwise,
   whence $\qu\kappa=\kappa$.
\qed

\abstand
\noindent
In the sequel let $\Kcal$ be the subset of \PG53 described in Theorem
\ref{theo-a}.
   \begin{remark}
   {\em
   By Theorem \ref{theo-d}, any set of points in \PG53 satisfying
   (\ref{condi}) and (\ref{condii}) is projectively equivalent to $\Kcal$. We
   infer from Lemma \ref{lemma-a} and Theorem \ref{theo-d} that the
   $12$--sets of points discussed in \cite{coxe58} and \cite{pell74} are
   essentially our $\Kcal$. By \cite[Teorema 4.3]{pell74}, conditions
   (\ref{condii}) and (\ref{condv}) characterize $\Kcal$ to within projective
   collineations. The set $\Kcal$ has a lot of fascinating geometric
   properties \cite{coxe58}, \cite{pell74}, \cite{todd59}.
   }
   \end{remark}

   \begin{remark}
   {\em
   Define a {\em block}\/ of $\Kcal$ as a hyperplane section of $\Kcal$
   containing more than three points. If $\Bcal$ denotes the set of all such
   blocks, then the incidence structure $(\Kcal,\Bcal,\in)$ is {\em Witt's\/
   $5$--$(12,6,1)$ design} $W_{12}$; cf., e.g., \cite[Chapter
   4]{beth-jung-lenz85}. According to Lemma \ref{lemma-a}, Theorem
   \ref{theo-c}, and Theorem \ref{theo-d}, such a point model of $W_{12}$ in
   \PG53 is projectively unique.
   }
   \end{remark}

   \begin{remark}
   {\em
   The automorphism group of $W_{12}$ is the {\em Mathieu group} $M_{12}$, a
   sporadic simple group acting sharply $5$--transitive on $\Kcal$; cf.,
   e.g., \cite[Chapter 4]{beth-jung-lenz85}. Each automorphism of
   $(\Kcal,\Bcal,\in)$ extends to a unique automorphic collineation of
   $\Kcal$ \cite{coxe58}, \cite{pell74}. Theorem \ref{theo-d} includes a short

   coordinate--free proof of that result.
   }
   \end{remark}

   \begin{remark}
   {\em
   The successive derivations of $W_{12}$ are a $4$--(11,5,1) design, a
   $3$--(10,4,1) design (the {\em M\"obius plane}\/ over the field extension
   $\GF9/ \GF3$), and a $2$--(9,3,1) design (the {\em affine plane}\/ over
   \GF3). One may obtain point models for them by suitable projections of
   $\Kcal$. Projection through a point of $\Kcal$ yields an $11$--cap in a
   hyperplane of \PG53. See \cite{hirs91}, \cite{pell73}, \cite{pell74},
   \cite{tall61}. If the centre of projection is a bisecant of $\Kcal$, then
   one gets an {\em elliptic quadric}\/ in a solid of \PG53. Finally, if the
   centre of projection is spanned by a triangle of $\Kcal$, then an {\em
   affine subplane}\/ of a projective plane of \PG53 arises. If the triangle
   is chosen to be $\Delta_\infty$, then there exists an affinity of this
   affine plane onto $\Pcal(\Vbf)\setminus l_\infty$. This is immediately
   seen from (\ref{vero}) and (\ref{delta}) by projecting, e.g., onto the
   plane with equations $x_{11}=x_{12}=x_{22}=0$.
   }
   \end{remark}

   \begin{remark}
   {\em
   Let $F^{\Pcal(\Wbf)}$ be the $F$--vector space of all functions
   $\Pcal(\Wbf)\rightarrow F$. Given $\Mcal\subset \Pcal(\Wbf)$ denote by
   $\chi(\Mcal)\in F^{\Pcal(\Wbf)}$ its characteristic vector (function).
   With the notations of Theorem \ref{theo-a} we obtain
      \begin{displaymath}
      \chi(\im\varphi) - \chi(\Gamma_\infty) + \chi(\Delta_\infty) =
      \chi(\Kcal).
      \end{displaymath}
   The characteristic vectors of the hyperplanes $\Hcal\subset\Pcal(\Wbf)$
   are spanning a linear $[364,22,121]$--code \cite[Theorem
   5.7.1]{assm-key92}. By (\ref{condi}), $\chi(\Kcal)$ is a word of weight
   $12$ in the orthogonal (dual) code, where orthogonality is understood with
   respect to the standard dot product. According to (\ref{condiii}), the
   Veronese variety yields a word of weight $13$ which has dot product $1\in
   F$ with each hyperplane. Thus, in terms of characteristic vectors, $\Kcal$
   arises from the Veronese variety by adding a word of weight $7$ which has
   dot product $2\in F$ with each hyperplane.

   Next let $\wbf_1,\ldots,\wbf_{12}\in\Wbf$ be vectors representing the
   points of $\Kcal$. As $f$ ranges over the dual vector space $\Wbf^\ast$,
   the words $(\wbf_1^f,\ldots,\wbf_{12}^f)\in F^{12}$
   give the {\em extended ternary Golay code\/} $\Gbf_{12}$. Cf.\
   \cite{assm-matt66}, where the dual point of view has been adopted.
   If we start instead with vectors $\vbf_1\vee\vbf_1, \ldots,
   \vbf_{13}\vee\vbf_{13}$ ($\vbf_i\in\Vbf$) representing the points of the
   Veronese surface, then we obtain a ternary $[13,6,6]$--code $\Cbf$,
   as follows from $\spn\im\varphi = \Pcal(\Wbf)$ and (\ref{condiii}).

   Given $f\in\Wbf^\ast$ then $\Abb{q}{\Vbf}{F}{\abf}{(\abf\vee\abf)^f}$ is
   a quadratic form. The mapping $f\mapsto q$ is a linear bijection of
   $\Wbf^\ast$ onto the vector space of quadratic forms $\Vbf\rightarrow F$.
   Thus, as $q$ ranges over all quadratic forms on $\Vbf$, the words
   $(\vbf_1^q, \ldots, \vbf_{13}^q)$ too comprise the code $\Cbf$.

   In order to identify the code $\Cbf$, let $\Cbf(p)$ ($p$ prime) be the
   linear code over \GF{p} which is spanned by the characteristic vectors of
   the lines of \PG2p. The dimension of $\Cbf(p)$ is $(p^2+p+2)/2$,
   $\Cbf(p)^\perp\subset\Cbf(p)$, and $\Cbf(p)^\perp$ has codimension $1$ in
   $\Cbf(p)$ \cite[49]{assm-key92}. Moreover, $\Cbf(p)^\perp$ coincides with
   two other codes arising from \PG2p: One is the code $\Ebf(p)$ spanned by
   the differences of characteristic vectors of lines \cite[Theorem
   6.3.1]{assm-key92}, the other is the code $\Cbf'(p)$ spanned by the
   characteristic vectors of the complements of lines, as follows easily from
   $\Cbf'(p)\subset\Cbf(p)^\perp$ and $\dim\Cbf'(p) = \dim\Cbf(p)^\perp$
   \cite[366]{brow-wilb95}.

   If a quadratic form \abb{q}{\Vbf}{F} is applied to four vectors $\vbf_i$
   which represent the points of a line, then one of the following
   (unordered) quadruples arises: $(0,0,0,0)$, $\pm(1,1,1,0)$, $(1,2,0,0)$,
   $(1,1,2,2)$. This is immediate from \cite[Lemma 5.2.1]{hirs79}. Hence
   $\Cbf\subset\Cbf(3)^\perp$ and, by $\dim\Cbf = \dim\Cbf(3)^\perp$, the two
   codes turn out to be the same.

   So, the self--dual extended ternary Golay code $\Gbf_{12}=\Gbf_{12}^\perp$
   is closely related to a self--orthogonal code
   $\Cbf\subset\Cbf^\perp=\Cbf(3)$ which belongs to an infinite family of
   codes obtained from \PG2p.
   }
   \end{remark}

   \begin{remark}
   {\em
   We aim at representing the points of $\Delta_\infty$ on the line
   $l_\infty$ by making use of the Veronese mapping $\varphi$: Each bijection
   of $l_\infty$ is a projectivity. There are three {\em elliptic
   involutions}\/ on $l_\infty$, each interchanging the points of $l_\infty$
   in pairs. Transformation under $\varphi$ yields three elliptic involutions
   on the conic $\Gamma_\infty$. Each of them extends uniquely to a harmonic
   homology of the plane $\Ecal_\infty$ leaving $\Gamma_\infty$ fixed, as a
   set \cite[2.4.4]{brau-I76}. The centres of the three homologies are three
   distinct internal points of $\Gamma_\infty$, whence they comprise the set
   $\Delta_\infty$. Thus the points of $\Delta_\infty$ are in one--one
   correspondence with the three elliptic involutions on $l_\infty$.

   Now it is natural to ask for a description of $W_{12}$ in terms of
   the nine points in $\Pcal(\Vbf)\setminus l_\infty$ and the three elliptic
   involutions on $\l_\infty$. It turns out that one obtains L\"uneburg's
   description \cite[Chapter 7]{luen69}, although from a different point of
   view. A block is precisely one of the following:
      \begin{enumerate}
      \item An affine line plus all three elliptic involutions.
      \item An ellipse together with those two elliptic involutions which are
      {\em not}\/ the involution of conjugate points on $\l_\infty$ with
      respect to the ellipse.
      \item A union of two distinct parallel affine lines.
      \item A cross of affine lines together with that elliptic involution
      which interchanges the points at infinity of the two lines.
      \end{enumerate}
   Cf.\ the proof of Theorem \ref{theo-a}. Thus each block arises from an
   affine quadric and certain elliptic involutions which are affine
   invariants of the quadric. This observation was the starting point for the
   present paper.
   }
   \end{remark}

\abstand

\end{document}